\documentclass{article}

\usepackage{PRIMEarxiv}
\usepackage[utf8]{inputenc} 
\usepackage[T1]{fontenc}    
\usepackage{hyperref}       
\usepackage{url}            
\usepackage{booktabs}       
\usepackage{amsfonts}       
\usepackage{nicefrac}       
\usepackage{microtype}      
\usepackage{lipsum}
\usepackage{fancyhdr}       
\usepackage{graphicx}       

\usepackage{subfig}
\usepackage{empheq}
\usepackage{longtable}
\usepackage{amsmath}
\usepackage{amsmath,bm}
\usepackage{amsthm}
\usepackage{float}
\usepackage{tabularx}
\usepackage{booktabs}
\usepackage{url}
\usepackage{tikz}
\usepackage{subfig}
\usepackage{caption}
\usepackage{siunitx}
\usepackage{float}

\newcommand{\numDOFs}{\ensuremath{N_p}}              
\newcommand{\numSnaps}{\ensuremath{N_s}}             
\newcommand{\numReducedBasis}{\ensuremath{L}}    

\newcommand{\param}{\ensuremath{\vectG{\zeta}}}      
\newcommand{\reducedCoeffV}{\ensuremath{\vect{\mathsf{y}}}}   
\newcommand{\reducedBasisV}{\ensuremath{\vect{\mathsf{v}}}}   
\newcommand{\reducedBasisT}{\ensuremath{\tensL{V}}}   

\newcommand{\vect}[1]{\ensuremath{\bm{#1}}}
\newcommand{\vectG}[1]{\ensuremath{\boldsymbol{#1}}} 	
\newcommand{\tensL}[1]{\ensuremath{\bm{\mathsf{#1}}}}   
\newcommand{\tensG}[1]{\ensuremath{\boldsymbol{\mathsf{#1}}}}	
\newcommand{\x}{\ensuremath{\vect{x}}}       
\newcommand{\normal}{\ensuremath{\vect{n}}} 
\newcommand{\room}{\ensuremath{\mathbb{R}}}
\newcommand{\domainDB}{\ensuremath{\Omega_D}}   
\newcommand{\domainRB}{\ensuremath{\Omega_R}}   
\newcommand{\solution}{\ensuremath{s}}   
\newcommand{\snapshotV}{\ensuremath{\vect{s}}}       
\newcommand{\snapshotT}{\ensuremath{\tensL{S}}}     
\newcommand{\reducedSnapshot}{\ensuremath{\vect{s}_L}}       
\newcommand{\interpolationFctFE}{\ensuremath{\phi}}   
\newcommand{\leftSVT}{\ensuremath{\tensL{W}}}   
\newcommand{\rightSVT}{\ensuremath{\tensL{Z}}}   
\newcommand{\EigenValT}{\ensuremath{\tensG{\Sigma}}}   
\newcommand{\squaredEigenValT}{\ensuremath{\tensG{\Lambda}}}   

\newcommand{\thermCond}{\ensuremath{\lambda}}          
\newcommand{\heatTransferCoef}{\ensuremath{\kappa}}          
\newcommand{\specHeat}{\ensuremath{c}}          
\newcommand{\temperature}{\ensuremath{T}}          
\newcommand{\temperatureAmbient}{\ensuremath{T_A}}          
\newcommand{\density}{\ensuremath{\rho}}   
\newcommand{\advectionVel}{\ensuremath{\vect{a}}}  
\newcommand{\grad}{\ensuremath{\vectG{\nabla}}}   
\newcommand{\domain}{\ensuremath{\Omega}}   

\newcommand{\inputnum}{1} 
\newcommand{\hiddennum}{5}  
\newcommand{\hiddenlayers}{2}  
\newcommand{\outputnum}{3} 
\newcommand{\colorInput}{white}
\newcommand{\colorOutput}{white}
\newcommand{\colorHidden}{white}
\newcommand{\colorLine}{black}

\title{A DATA-DRIVEN REDUCED ORDER MODELING APPROACH APPLIED IN CONTEXT OF NUMERICAL ANALYSIS AND OPTIMIZATION OF PLASTIC PROFILE EXTRUSION}

\author{
	Daniel Hilger\\
	Chair for Computational Analysis of Technical Systems\\
	RWTH Aachen University\\
	Schinkelstr. 2, 52062, Aachen\\
	\texttt{hilger@cats.rwth-aachen.de} \\
	\textit{www.cats.rwth-aachen.de}
	\AND
	Norbert Hosters\\
	Chair for Computational Analysis of Technical Systems\\
	RWTH Aachen University\\
	Schinkelstr. 2, 52062, Aachen\\
	\texttt{hosters@cats.rwth-aachen.de} \\
	\textit{www.cats.rwth-aachen.de}
	
}

\theoremstyle{remark}

\begin{document}
	\maketitle
\begin{abstract}
In course of this work, we examine the process of plastic profile extrusion, where a polymer melt is 
shaped inside the so-called extrusion die and fixed in its shape by solidification in the downstream 
calibration unit.
More precise, we focus on the development of a data-driven reduced order model (ROM) for the purpose of predicting temperature distributions within the extruded profiles inside the calibration unit.
Therein, the ROM functions as a first step to our overall goal of prediction based process control in order to avoid undesired warpage and damages of the final product. 
\end{abstract}
\section{INTRODUCTION}
\begin{figure}[h!]
	\centering
	\input{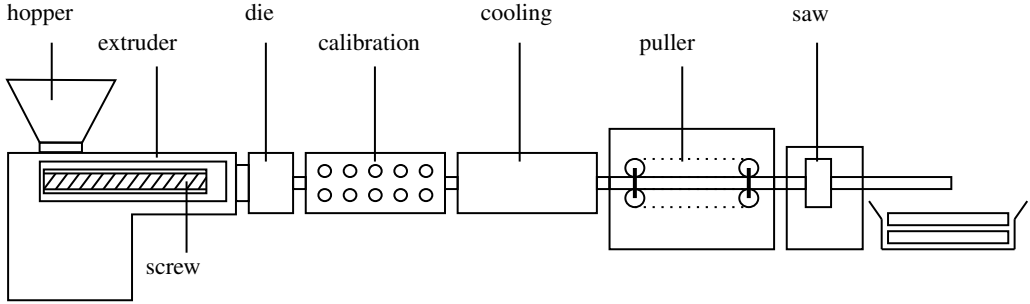}
	\caption{Process sections in classic extrusion process.}
	\label{fig:SketchExtrusionLine}
\end{figure}

Nowadays, production processes are highly optimized, monitored and controlled.
Optimization and control need a lot of information about the process which can partially be provided by sensor data.
Measurements of these sensors, however, can only be made at selected points or from outside of the product.
For this reason and because the overall number of sensors that can be used within a process is generally limited, the knowledge about the process is always incomplete.
One way to complete these knowledge gaps can be overcome by additional information provided by numerical analysis.
Up-to-date, numerical analysis is already widely used in process design and optimization.
The computational methods applied in numerical analysis, as the finite element method, finite volumes or finite differences, just to name a few, are already well established and provide good predictions about the processed materials.
These methods, however, have all in common that they become quickly costly in terms of computational time, even when solved on so-called high performance computing (HPC) architectures.
Time-consuming simulations cannot be used in process control since here time is a limiting factor.
Therefore, the development of reduced order models (ROM) has become of particular interest over the last decades.\\
\\
There exist already a variety of ROM approache, which in literature are often categorized into hierarchical, projection based or data-driven ROM \cite{EldredDunlavy2006,BennerGugercinWillcox2015}.
In hierarchical ROM, simply the physics of the problem is reduced, e.g. by considering only the main flow direction \cite{FringsBehrElgeti2017}.
In projection based methods the high dimensional solution space is projected onto a smaller, reduced space.
The problem is then solved on the reduced solution space.
Projection based methods require access to the source code and often need a  problem dependent implementation.
Examples to this methodology can be found in \cite{HesthavenRozzaStamm2016,KeyEtAl2021}.
In scope of this work, we will though focus on the remaining group, namely the data-driven approach.\\
\\
The category of data-driven ROM approaches, once again, incorporates a variety of different strategies.
Especially increased interest in machine learning evoked a larger number of new data-driven ROM strategies as for example physics informed neural networks (PINN) \cite{RaissiPerdikarisKarniadakis2019} or convolutional autoencoder \cite{FukamiTaichiKoji2022}.
The strategy we pursue within the presented work is based on the ROM strategy proposed by Hesthaven and Ubbiali \cite{HesthavenUbbiali2018}.
Therein, the solution to a problem formulation is approximated under consideration of existing data sets, where each data set is associated with a distinct set of process parameters.
Based on this data, we apply the so-called proper orthogonal decomposition (POD) \cite{Chatterjee2000} to identify a reduced basis (RB) representation for the data set.
Once the  RB is identified, it is possible to recover reduced coefficients for a specific parameter setting.
The relation between process parameters and reduced coefficient is then what is trained in the ROM to approximate the solution of an unseen set of process parameters.
Possible strategies to the create this mapping are radial basis functions (RBF) \cite{WaltonHassanMorgan2013,XiaoEtAl2015}, Gaussian progress regression (GPR)\cite{RozzaEtAl2018} or artificial neural networks (ANN) \cite{HesthavenUbbiali2018,BerzinsEtAl2020,WangHesthavenRay2019}.\\
\\
In the presented work, we apply a RB-ANN ROM to the industrial process of plastic profile extrusion.
Inside extrusion process lines, as sketched in \ref{fig:SketchExtrusionLine}, granular, raw plastics are continuously processed to profiles with fixed cross-sectional shape.
The extrusion process can therein be divided into two section, first the hot mixing and shaping part and the subsequent cooling and calibration part.
In the remaining work, we will focus on the cooling and calibration part of extrusion.
Inside the calibration unit the still liquid plastic mold is cooled down so that it solidifies and is fixed in shape.
Thereby, the cooling of the extrudate should be as uniform as possible, since otherwise undesired warpage and distortion will damage the profile.
The ROM presented in this work, should be cable of delivering advance information on the temperature distribution inside the profile and
thereby function as decision support tool in process adaptation under varying process conditions.

\section{MATERIAL AND NUMERICAL MODELING}\label{sec:problem}
To fully describe the behavior of plastics melts in the calibration unit, it would require the conservation equation for mass, momentum and energy, as well as further constitutive models to describe the non-Newtonian nature of plastics, depending on the type of plastic also crystallization models and finally a methodology that is able to represent the domain changes due to warpage and deformation. 
This work, however, wants to focus on the evaluation and exploration of the ROM aspect.
Therefore, we consider only the temperature equation to describe the cooling plastic melt inside.
The temperature distribution $\temperature$ within the calibration unit, represented by $\domain$, is described by the following equation:
\begin{align}
	\density \specHeat \advectionVel \grad \temperature - \thermCond \grad^2 \temperature = 0\qquad&\text{in}~\domain. \label{eq:Heat}
\end{align}
Here, $\density$ represents the density, $\specHeat$ the specific heat, $\advectionVel$ the advection velocity and $\thermCond$ the thermal conductivity.
Uniform temperatures can be assumed as Dirichlet conditions at the inlet of the calibriation unit, whereas the emitted heat flux through the profile surface depends on the difference between ambient temperature $\temperatureAmbient$ and surface temperature.
Further, the heat transfer coefficient between the processed material and the surrounding cooling fluid influences the heat flux through the surface:
\begin{align}
	\temperature = \temperature_{in}, &\qquad on~\domainDB,\\
	\normal \cdot \heatTransferCoef \grad \temperature = \heatTransferCoef \left(\temperature - \temperatureAmbient\right) , &\qquad on~\domainRB.\label{eq:RBC}
\end{align}

The data used within our ROM approach is generated with our in-house HPC stabilized finite element code.
For a more detailed discussion about utilized discretization and stabilization methods, we refer the interested reader to \cite{HelmigBehrElgeti2019}.

\section{STANDARDIZED NON-INTRUSIVE REDUCED ORDER MODEL}\label{sec:method}
The ROM used within this work is the so-called standardized non-intrusive reduced order model (sniROM) propsed by B\={e}rzi\c{n}\v{s} et al. \cite{BerzinsEtAl2020}.
It follows the offline-online paradigm presented by Ubbiali and Hesthaven \cite{HesthavenUbbiali2018}.
Within the offline step, the underlying data for the ROM is aggregated first.
Subsequently, this data is used to construct a RB.
Lastly, the data and the RB are utilized to train an interpolation model, which establishes a relationship between problem specific parameters and a linear combination of RB vectors.
In the online stage, the ROM can then  be evaluated quickly for unseen parameters at low computational cost.

\subsection{Solution representation via reduced basis}
The pre-computed data consist of pointwise stored solutions to equations \eqref{eq:Heat}-\eqref{eq:RBC} for varying parameter sets $\param$.
In scope of this work, $\param$ is a vector containing the ambient temperature and the heat transfer coefficient  $\heatTransferCoef$.
The pointwise stored solutions will be referred to synonymously as snapshots in the remaining work.
Similar as in the finite element method (FEM), the pointwise solution $\solution$ can be represented on the domain by means of  spatial basis functions $\interpolationFctFE$:
\begin{align}
	\solution \left(\param;\x\right) = \snapshotV\left(\param\right)^T\interpolationFctFE\left(\x\right) = \sum_{i=1}^{\numDOFs} \snapshotV_i\left(\param\right)\interpolationFctFE_i\left(\x\right),
\end{align}
where $\numDOFs$ denotes the number of data points per parameter configuration.
The principle idea of reduced basis methods (RBM) is now to approximate any solution vector $\snapshotV \left(\param\right)$ by a linear combination of the first $\numReducedBasis$ RB vectors $\reducedBasisT = \left[\reducedBasisV^1 | \dots | \reducedBasisV^\numReducedBasis\right]$ of a given dataset:
\begin{align}
	\snapshotV_i\left(\param\right) \approx \reducedSnapshot\left(\param\right) =  \sum_{i=1}^{\numReducedBasis} \reducedCoeffV_l \reducedBasisV^l = \reducedBasisT \reducedCoeffV \left(\param\right),
\end{align}
where $\reducedCoeffV$ denote the reduced coefficients.
With the pre-computed snapshots, the reduced coefficients to a corresponding parameterset can be recovered via the RB vectors $\reducedBasisT$.
\begin{equation}
	\reducedCoeffV\left(\param\right) = \reducedBasisT^T \snapshotV\left(\param\right)
\end{equation}

\subsection{Construction of reduced basis}
\begin{figure}[h!]
	\begin{tikzpicture}[scale=1.0]
		\node (0,0) {\includegraphics[width=1.0\textwidth,trim=6cm 1cm 6cm 1cm,clip]{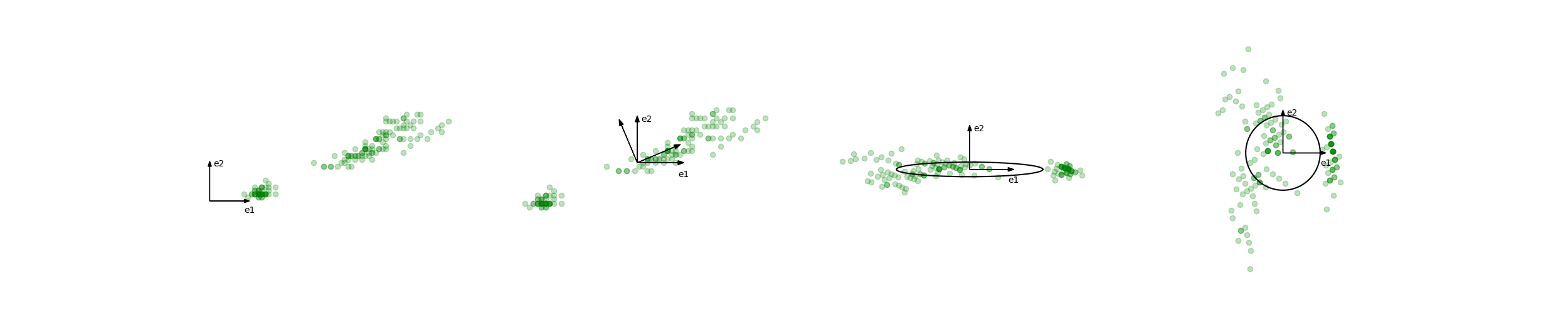}};
		\draw [-latex,line width=1pt](-4.5,-1.2)  -- (-3.5,-1.2) node[midway,below]{zero-center};
		\draw [-latex,line width=1pt](-0.5,-1.2) -- (0.5,-1.2)node[midway,below]{decorrelate};
		\draw [-latex,line width=1pt](4.0,-1.2) -- (5.0,-1.2)node[midway,below]{normalize scale};
	\end{tikzpicture}
	\caption{Data preprocessing.}
	\label{fig:DataStandardization}
\end{figure}
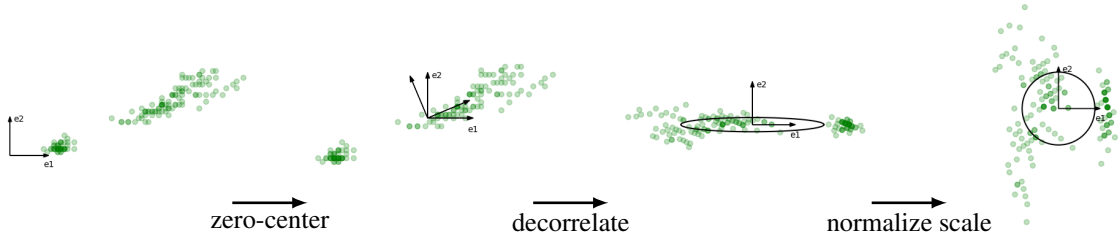

The RB is constructed from a pre-computed dataset.
The dataset consist of snapshots of a distinct parameter set.
These snapshots form the columns to the overall snapshot matrix $\snapshotT \in \room^{\numDOFs \times \numSnaps}$, where $\numSnaps$ is the total number of snapshots used to construct the RB.
Performing a proper orthogonal decomposition (POD) \cite{Chatterjee2000} on the normalized snapshot matrix yields:
\begin{align}
	\snapshotT / \sqrt{\numSnaps} = \leftSVT \EigenValT \rightSVT^T,
	\label{eq:reducedCoefficients}
\end{align}
where, $\EigenValT$ is a diagonal matrix containing the descending eigenvalues of $\snapshotT$, whereas $\leftSVT$ and $\rightSVT$ represent the right and left singular vectors.
The first $\reducedBasisV$ left singular Eigenvectors $\leftSVT$ are chosen to be the RB of our ROM.
Following the suggestion in \cite{BerzinsEtAl2020} the snapshot data should be pre-processed, cmp. Figure \ref{fig:DataStandardization}.
The data is centered and parameters and coefficients are standardized.
However, the POD is not the most efficient strategy to compute the RB.
It can be computed more efficiently by the so-called method of snapshots \cite{Sirovich1991},
where the dimensionality of our dataset is exploited. 

\paragraph{Method of snapshots}
The datasets based on FE simulations, typically has significantly less snapshots than degrees of freedom per snapshot.
Since the squared snapshot matrix $\snapshotT\snapshotT^T/{\numSnaps} $ is only of dimension $\room^{\numSnaps \times \numSnaps}$, it is computational less expensive to compute the reduced basis of $\snapshotT$ via singular value decomposition of its square:
\begin{align}
	\snapshotT\snapshotT^T/{\numSnaps}  = \rightSVT \squaredEigenValT \rightSVT^T.
\end{align}
This can then be used to calculate the eigenvalues as $\EigenValT =\squaredEigenValT^{1/2}$ and the left eigenvectors of $\snapshotT$ by a simple matrix product $\leftSVT = \snapshotT \rightSVT \EigenValT^{-1}/{\numSnaps}$

\subsection{Interpolation of reduced coefficients}
\begin{figure}[h!]
	\centering
	\resizebox{0.45\textwidth}{!}{
		\begin{tikzpicture}[scale=1]
			
			\node[circle,minimum size = 6mm, fill=\colorInput!70,draw] (Input-1) at (0,-1) {$\boldsymbol{\zeta}$};
			
			\node[circle, minimum size = 6mm, fill=\colorOutput!70,	yshift=(\outputnum-\inputnum)*5 mm,draw] (Output-1) at (2.5*\hiddenlayers+2.5,-1) {$y_1$};
			\node[circle, minimum size = 6mm, fill=\colorOutput!70,	yshift=(\outputnum-\inputnum)*5 mm,fill=none ] (Output-2) at (2.5*\hiddenlayers+2.5,-2) {\small $\begin{array}{c} \cdot\\ \cdot\\ \cdot\end{array}$};
			\node[circle, minimum size = 6mm, fill=\colorOutput!70,	yshift=(\outputnum-\inputnum)*5 mm,draw] (Output-3) at (2.5*\hiddenlayers+2.5,-3) {$y_L$}; 		
			
			\foreach \j in {1,...,\hiddenlayers}
			{
				\foreach \i in {1,...,\hiddennum}
				{
					\node[circle, 
					minimum size = 6mm,
					fill=\colorHidden!50,
					yshift=(\hiddennum-\inputnum)*5 mm,
					draw
					] (Hidden-\j-\i) at (2.5*\j,-\i) {};
				}
			}
			\foreach \i in {1,...,\inputnum}
			{
				\foreach \j in {1,...,\hiddennum}
				{
					\draw[\colorLine,->, shorten >=1pt] (Input-\i) -- (Hidden-1-\j);   
				}
			}
			
			\foreach \i in {1,...,\hiddennum}
			{
				\foreach \j in {1,...,\hiddennum}
				{
					\draw[\colorLine,->, shorten >=1pt] (Hidden-1-\i) -- (Hidden-2-\j);   
				}
			}
			\foreach \i in {1,...,\hiddennum}
			{
				\foreach \j in {1,...,\outputnum}
				{
					\draw[\colorLine,->, shorten >=1pt] (Hidden-2-\i) -- (Output-\j);
				}
			}
			
	\end{tikzpicture}}
	
	\caption{Examplary NN used to reduced coefficient interpolation.}
	\label{fig:NN}
\end{figure}
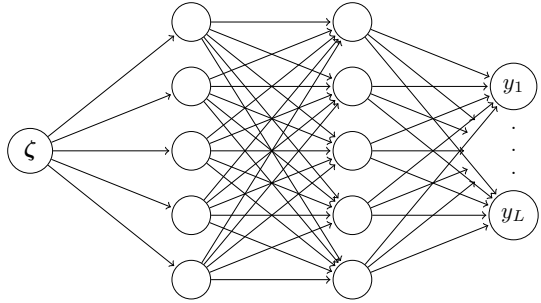

Once the reduced basis vectors are computed to a data set, equation \eqref{eq:reducedCoefficients} is used to identify the reduced coefficients for each snapshot vector in $\snapshotT$. 
This mapping from distinct parameters to a reduced coefficients vector is then utilized to train an interpolation model so that the interpolation model proposes reduced coefficients for an unseen choice of process parameters.
The interpolation model chosen within this work is a feedforward neural network (FNN), as shown in Figure \ref{fig:NN}.
The inputs to our network are the varying process parameters $\param$, whereas the reduced coefficients $\reducedCoeffV$ are the ouputs.
The number of hidden layers and their size are part of the so-called hyperparameters.
The FNN is implemented with hyperbolic tangent activation functions.
All parameters that describe network or specify the training configuration are called hyperparameters.
These parameters affect the quality of the interpolation model.
Finding the combination of parameters that result the best interpolation is part of the offline training process.

\section{TEMPERATURE DISTRIBUTION IN CALIBRATION UNIT FOR EXTRUDED 3D-DUMBBELL PROFILE}
\begin{figure}[H]
	\centering
	\resizebox{\textwidth}{!}{
		\tikzset{every picture/.style={line width=0.75pt}} 
		\begin{tikzpicture}[x=0.75pt,y=0.75pt,yscale=-1,xscale=1]
			
			\draw  [draw opacity=0][line width=1.5]  (159.66,160.04) .. controls (154.24,175.7) and (139.36,186.94) .. (121.86,186.94) .. controls (99.77,186.94) and (81.86,169.03) .. (81.86,146.94) .. controls (81.86,124.84) and (99.77,106.94) .. (121.86,106.94) .. controls (138.02,106.94) and (151.95,116.52) .. (158.26,130.32) -- (121.86,146.94) -- cycle ; \draw  [line width=1.5]  (159.66,160.04) .. controls (154.24,175.7) and (139.36,186.94) .. (121.86,186.94) .. controls (99.77,186.94) and (81.86,169.03) .. (81.86,146.94) .. controls (81.86,124.84) and (99.77,106.94) .. (121.86,106.94) .. controls (138.02,106.94) and (151.95,116.52) .. (158.26,130.32) ;  
			\draw [line width=1.5]    (158,130) -- (238,130) ;
			\draw [line width=1.5]    (159.66,160.04) -- (238,160) ;
			\draw  [draw opacity=0][line width=1.5]  (238.14,130.03) .. controls (243.36,121.59) and (252.64,115.93) .. (263.29,115.8) .. controls (279.86,115.59) and (293.46,128.85) .. (293.66,145.42) .. controls (293.87,161.99) and (280.6,175.59) .. (264.04,175.79) .. controls (252.59,175.93) and (242.55,169.64) .. (237.38,160.27) -- (263.67,145.79) -- cycle ; \draw  [line width=1.5]  (238.14,130.03) .. controls (243.36,121.59) and (252.64,115.93) .. (263.29,115.8) .. controls (279.86,115.59) and (293.46,128.85) .. (293.66,145.42) .. controls (293.87,161.99) and (280.6,175.59) .. (264.04,175.79) .. controls (252.59,175.93) and (242.55,169.64) .. (237.38,160.27) ;  
			\draw [line width=0.75]    (123,145) -- (81.8,145.2) ;
			\draw [line width=0.75]    (294.2,145.6) -- (263.67,145.79) ;
			\draw [line width=0.75]    (238,180) -- (158,180) ;
			\draw [line width=0.75]    (158,175) -- (158,185) ;
			\draw [line width=0.75]    (238,175) -- (238,185) ;
			\draw [line width=0.75]    (148,160) -- (148,130) ;
			\draw [line width=0.75]    (153,130) -- (143,130) ;
			\draw [line width=0.75]    (153,160) -- (143,160) ;
			\draw  [dash pattern={on 0.84pt off 2.51pt}]  (390,115) -- (600,115) ;
			\draw  [dash pattern={on 0.84pt off 2.51pt}]  (390,175) -- (600,175) ;
			\draw  [dash pattern={on 0.84pt off 2.51pt}]  (390,160) -- (600,160) ;
			\draw  [dash pattern={on 0.84pt off 2.51pt}]  (390,130) -- (600,130) ;
			\draw [line width=0.75]    (600,200) -- (390,200) ;
			\draw [line width=0.75]    (390,195) -- (390,205) ;
			\draw [line width=0.75]    (600,195) -- (600,205) ;
			\draw [line width=1.5]    (389,187) -- (600,187) -- (600,103) -- (389,103) -- cycle ;
			\draw    (50,148) -- (50,175) ;
			\draw [shift={(50,145)}, rotate = 90] [fill={rgb, 255:red, 0; green, 0; blue, 0 }  ][line width=0.08]  [draw opacity=0] (8.93,-4.29) -- (0,0) -- (8.93,4.29) -- cycle    ;
			\draw    (77,175) -- (50,175) ;
			\draw [shift={(80,175)}, rotate = 180] [fill={rgb, 255:red, 0; green, 0; blue, 0 }  ][line width=0.08]  [draw opacity=0] (8.93,-4.29) -- (0,0) -- (8.93,4.29) -- cycle    ;
			\draw    (351,148) -- (351,175) ;
			\draw [shift={(351,145)}, rotate = 90] [fill={rgb, 255:red, 0; green, 0; blue, 0 }  ][line width=0.08]  [draw opacity=0] (8.93,-4.29) -- (0,0) -- (8.93,4.29) -- cycle    ;
			\draw    (378,175) -- (351,175) ;
			\draw [shift={(381,175)}, rotate = 180] [fill={rgb, 255:red, 0; green, 0; blue, 0 }  ][line width=0.08]  [draw opacity=0] (8.93,-4.29) -- (0,0) -- (8.93,4.29) -- cycle    ;
			
			\draw (91,124) node [anchor=north west][inner sep=0.75pt]   [align=left] {$r_1$};
			\draw (267,124) node [anchor=north west][inner sep=0.75pt]   [align=left] {$r_2$};
			\draw (152,134) node [anchor=north west][inner sep=0.75pt]   [align=left] {$h$};
			\draw (192,184) node [anchor=north west][inner sep=0.75pt]   [align=left] {$w$};
			\draw (484,205) node [anchor=north west][inner sep=0.75pt]   [align=left] {$l$};
			\draw (59,180) node [anchor=north west][inner sep=0.75pt]   [align=left] {$x$};
			\draw (33,141) node [anchor=north west][inner sep=0.75pt]   [align=left] {$y$};
			\draw (360,180) node [anchor=north west][inner sep=0.75pt]   [align=left] {$z$};
			\draw (334,141) node [anchor=north west][inner sep=0.75pt]   [align=left] {$y$};
			
		\end{tikzpicture}
	}
	\caption{Sketch of dumbbell profile.}
	\label{fig:Sketch}
\end{figure}
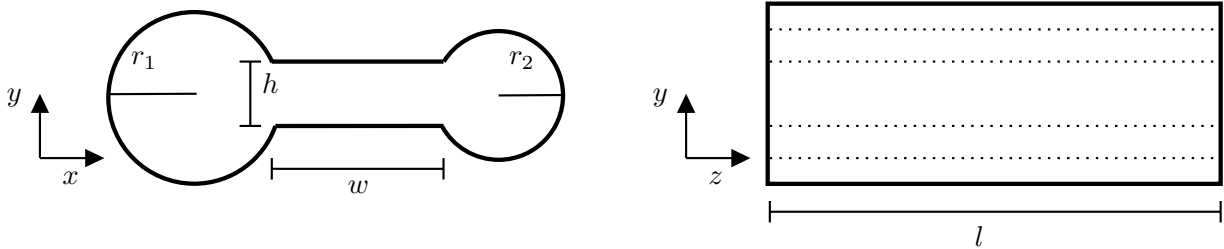

In this section we apply, as an example, the model presented in section \ref{sec:method}, for the problem described in section \ref{sec:problem}, to a so-called dumbbell extrusion profile.
The dimensions of the utilized geometry are given in Figure \ref{fig:Sketch} and in Table \ref{tab:dimensions}.
The material parameters used within the FE simulations for the data aggregation are shown in Table \ref{tab:processParameters}.
The parameters varied in the data set are the ambient temperature $\temperatureAmbient$ and the heat transfer coefficient  $\heatTransferCoef$ between material and surrounding fluid.

\begin{table}[H]
	\centering
	\begin{tabular}{llll}
		\multicolumn{3}{l}{\textbf{Geometric dimensions}}  \\
		\hline
		$r_1 $& radius  & $0.006 $ &   $\left[\mathrm{m}\right]$   \\
		$r_2$ & radius  & $0.0033   $ &   $\left[\mathrm{m}\right]$      \\
		$h$    & height  &$ 0.01   $   &   $\left[\mathrm{m}\right]$     \\
		$w$    & width   &$ 0.003    $&   $\left[\mathrm{m}\right]$       \\ 
		$l$    & length   &$ 1.0    $ &   $\left[\mathrm{m}\right]$  
	\end{tabular}
	\caption{Geometric parameters.}
	\label{tab:dimensions}
\end{table}
\begin{table}[H]
	\center
	\begin{tabular}{lllr}
		\multicolumn{4}{l}{\textbf{Process and material parameters}}  \\
		\hline
		$\density$ & density  & 900                    &   $\left[\mathrm{kg\,m^{-3}}\right]$                 \\
		$\advectionVel$ & advection velocity  & 0.00011     &   $\left[\mathrm{m\,s^{-1}}\right]$                               \\
		$\temperatureAmbient$    & ambient temperature  & $\left[288,298\right] $&   $\left[\mathrm{K}\right]$    \\
		$\heatTransferCoef$    & heat transfer coefficient   & $\left[-320,-218\right]$     &   $\left[\mathrm{W}/\left(\mathrm{m}^2\,\mathrm{K}\right)\right]$ \\
		& & & \\
		\multicolumn{4}{l}{\textbf{ROM and hyperparameters}}  \\
		\hline
		\multicolumn{3}{l}{size trainingsset} & $100$                    \\
		\multicolumn{3}{l}{size testset} & $100 $                   \\
		\multicolumn{3}{l}{size validationset} & $100 $              \\
		\multicolumn{3}{l}{number outputs/reduced basis vectors }& $30$                    \\
		\multicolumn{3}{l}{number hidden layers} & 10                    \\
		\multicolumn{3}{l}{number of neurons }& 40                    		                             
	\end{tabular}
	\caption{Process and model parameters.}
	\label{tab:processParameters}
\end{table}

\begin{figure}[h!]
	\centering
	\subfloat[Constant inlet BC. \label{fig:a}]{\includegraphics[width=0.4\textwidth]{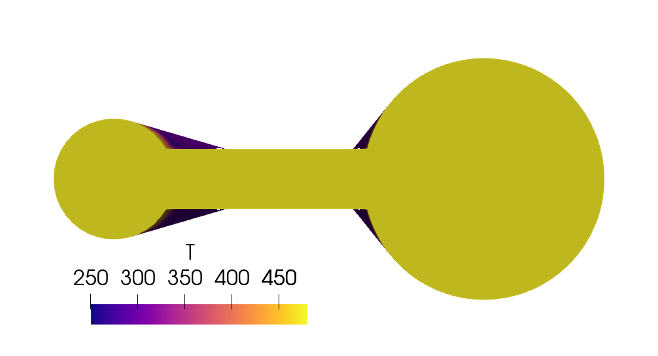}}\qquad
	\subfloat[Temperature prediction outlet \label{fig:a}]{\includegraphics[width=0.4\textwidth]{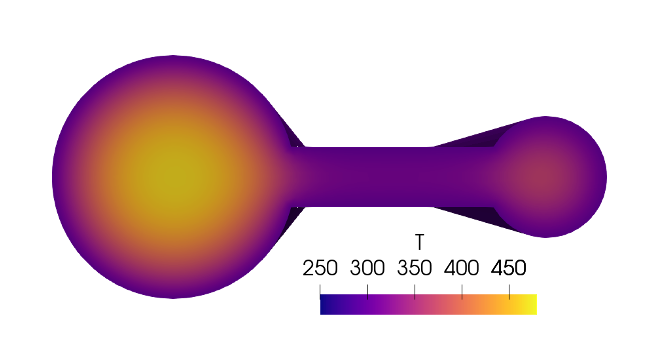}}
	\caption{Cross sectional temperature distributions at inlet (a)  and outlet (b))}
	\label{fig:InletOutlet}
\end{figure}

\begin{figure}[h!]
	\center
	\includegraphics[width=0.9\textwidth]{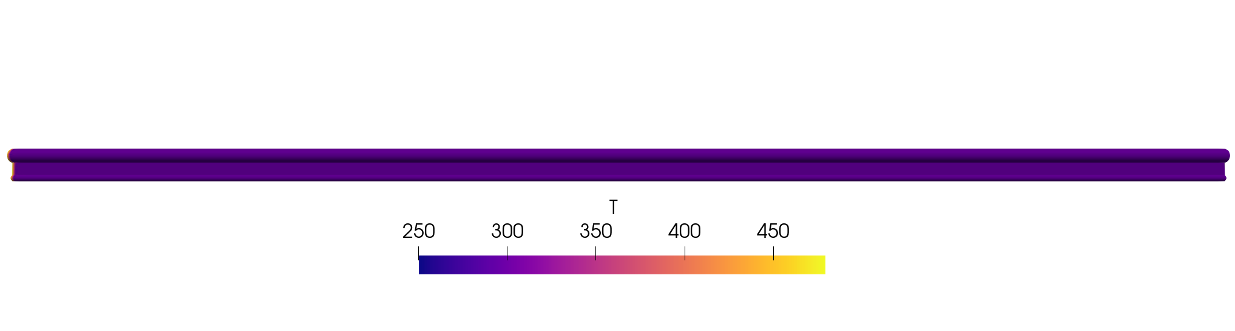}
	\caption{Predicted temperature on surface of extrudate.}
	\label{fig:tempDistributionSurface}
\end{figure}

\begin{figure}[h!]
	\center
	\includegraphics[width=0.9\textwidth]{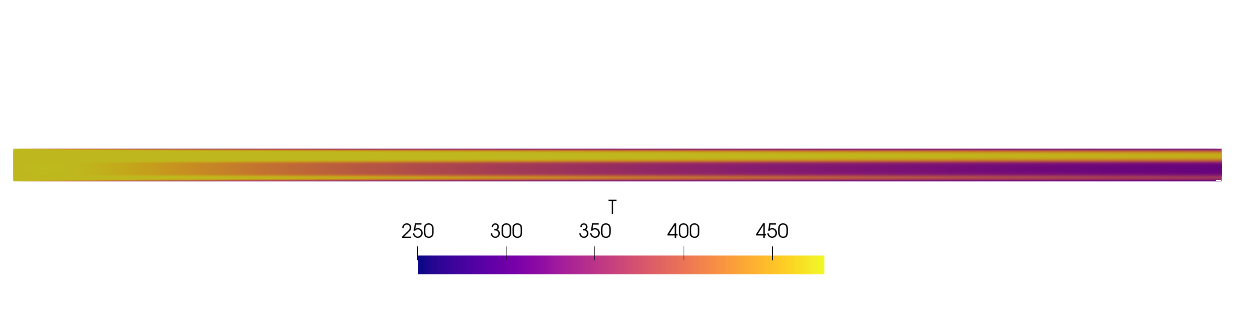}
	\caption{Predicted temperature inside extrudate.}
	\label{fig:tempDistributionCrossSection}
\end{figure}

The NN used within the ROM is trained with the hyperparameters listed in table \ref{tab:processParameters}.
The ROM temperature predictions on the dumbbell geometry, cmp. Figures \ref{fig:InletOutlet}-\ref{fig:tempDistributionCrossSection}, match the expected temperature distributions.
In Figures \ref{fig:InletOutlet} and \ref{fig:tempDistributionCrossSection}, it can be observed that the extruded plastic is still hotter in the larger cylinder of the dumbbell at the outlet, even though the surface temperature of the profile, cmp. Figure \ref{fig:tempDistributionSurface}, has cooled down to a uniform  temperature.
These results indicate that implemented ROM could be applied in process control with only small modifications, e.g. extending the model by an additional heat transfer coefficient to handle both sides of the profile differently.

\section{CONCLUSION}
In scope of this work, we presented a data-driven ROM approach.
The ROM allows to predict temperature distribution on plastic extrusion profiles for varying ambient temperatures and heat-transfer coefficients.
For the presented dumbbell example, we could show that the model predictions for unseen process configurations resulted physically reasonable results.
This work is only a first conceptual work and needs to extended by further investigations of the model.
ln future work, the relation between accuracy and size of trainingsset, as well as the model behavior for more complex geometries should be investigated.
Nevertheless, we can conclude that the model appears to be suitable for the application in process controlling, where a ROM can be created during the process desgin, similar as the extrusion tool itself.
\section*{ACKNOWLEDGMENT}
Funded by the Deutsche Forschungsgemeinschaft (DFG, German Research Foundation) under Germany’s Excellence Strategy–EXC-2023 Internet of Production–390621612.
Further, the authors gratefully acknowledge the computing time granted by the JARA Vergabegremium and provided on the JARA Partition part of the supercomputer JURECA at Forschungszentrum Jülich.

\bibliographystyle{ieeetr}  
\bibliography{Literature}
	
\end{document}